\documentstyle[12pt,draft]{article}

 \newcounter{abceqn}

 \newcounter{abcfig}

\newcommand{\na}{\nabla}

\newcommand{\Om}{\Omega}

\newcommand{\pa}{\partial}

\newcommand{\Dl}{\Delta}

\newcommand{\sg}{\sigma}

\newcommand{\La}{\Lambda}
\newcommand{\la}{\lambda}

\renewcommand{\thesection}{\arabic{section}.}

\newcommand{\eqnsection}[1]{
	\section{#1}
	\setcounter{equation}{0}
	\renewcommand{\theequation}{\thesection\arabic{equation}}
	\setcounter{figure}{0}
	\renewcommand{\thefigure}{\arabic{figure}}
	\setcounter{remark}{0}
	\renewcommand{\theremark}{\thesection\arabic{remark}}
	\setcounter{theorem}{0}
	\renewcommand{\thetheorem}{\thesection\arabic{theorem}}
	\setcounter{lemma}{0}
	\renewcommand{\thelemma}{\thesection\arabic{lemma}}
}

\title{{\bf Isospectral Theory of Euler Equations}}

\author{ \\ \\ \\ \\ 
Yanguang (Charles)\ \ Li  
\\  \\ Department of Mathematics, 
 \\ University of Missouri \\ 
Columbia, MO 65211 \\ E-mail: cli@math.missouri.edu \\ \\ \\
Roman Shvidkoy \\ \\ Department of Mathematics, 
 \\ University of Missouri \\ 
Columbia, MO 65211 \\ E-mail: shvidkoy@math.missouri.edu}

\date{\today}

\begin{document}
\bibliographystyle{unsrt}
\maketitle
\newpage
\begin{abstract}    
Isospectral problem of both 2D and 3D Euler equations of inviscid 
fluids, is investigated. Connections with the Clay problem are 
described. Spectral theorem of the Lax pair is studied. \\



\end{abstract}

\newtheorem{lemma}{Lemma}
\newtheorem{theorem}{Theorem}
\newtheorem{corollary}{Corollary}
\newtheorem{remark}{Remark}
\newtheorem{definition}{Definition}
\newtheorem{proposition}{Proposition}
\newtheorem{assumption}{Assumption}

\eqnsection{Introduction}

This note is a continuation of the works \cite{Li01a} and \cite{LY00}.
It focus upon the isospectral property of the Lax pairs of both 2D 
and 3D Euler equations. It provides efforts towards the connection 
between isospectral theory and the Clay problem on Navier-Stokes 
equations.

\eqnsection{2D Euler Equation}

The 2D Euler equation can be written in the vorticity form,
\begin{equation}
\pa_t \Om + \{ \Psi, \Om \} = 0 \ ,
\label{2deuler}
\end{equation}
where the bracket $\{\ ,\ \}$ is defined as
\[
\{ f, g\} = (\pa_x f) (\pa_y g) - (\pa_y f) (\pa_x g) \ ,
\]
where $\Psi$ is the stream function given by,
\[
u=- \pa_y \Psi \ ,\ \ \ v=\pa_x \Psi \ ,
\]
where $u$ and $v$ are respectively the velocity components 
along $x$ and $y$ directions, and the relation between 
vorticity $\Om$ and stream function $\Psi$ is,
\[
\Om =\pa_x v - \pa_y u =\Dl \Psi \ .
\]
\begin{theorem}[Li, \cite{Li01a}]
The Lax pair of the 2D Euler equation (\ref{2deuler}) is given as
\begin{equation}
\left \{ \begin{array}{l} 
L \varphi = \la \varphi \ ,
\\
\pa_t \varphi + A \varphi = 0 \ ,
\end{array} \right.
\label{laxpair}
\end{equation}
where
\[
L \varphi = \{ \Om , \varphi \}\ , \ \ \ A \varphi = \{ \Psi, \varphi \}\ ,
\]
and $\la$ is an imaginary constant, and $\varphi$ is a complex-valued function.
\label{2dlp}
\end{theorem}

\subsection{Isospectral Theory and Conservation Laws}

Denote by $H^s$ the Sobolev space $H^s(R^2)$ or $H^s(T^2)$, 
and $\| \ \|_s$ the $H^s$ norm.
\begin{theorem}
Let $\Om$ be a solution to the 2D Euler equation (\ref{2deuler}),
$\varphi$ be a solution to the Lax pair (\ref{laxpair}) at ($\Om , \la $),
then 
\[
I = \frac{\| \{ \Om , \varphi \} \|_s}{\| \varphi \|_s}
\]
is conservation law, i.e., I is independent of $t$. 
\end{theorem}

Proof: Take the $H^s$ norm on both side of the first equation 
in the Lax pair (\ref{laxpair}), then 
\[
I = \frac{\| \{ \Om , \varphi \} \|_s}{\| \varphi \|_s} = |\la |.
\]
By the isospectral property of the Lax pair (\ref{laxpair}), 
I is independent of $t$. $\Box$

An interesting idea is to try to use the conservation law $I$ to 
prove global well-posedness. In this 2D Euler case, the idea is not 
very important since the global well-posedness is already proved. 
The hope is to use this idea to prove the global well-posedness of
3D Euler equation.

\begin{lemma}
If $\varphi$ solves the Lax pair (\ref{laxpair}), then $f(\varphi)$ 
solves
\begin{equation}
\pa_t f(\varphi) + \{ \Psi , f(\varphi) \} = 0 \ ,
\label{tra}
\end{equation}
for any $f$ smooth enough.
\label{trale}
\end{lemma}

Proof: The proof is completed by direct verification. $\Box$

\begin{lemma}
If $\{ \phi_j \}_{j=1,2, \cdots }$ is a complete base of $H^s$,
where $\phi_j$'s are $f(\varphi)$'s at different values of $\la$;
then
\[
\Om = \sum_{j=1}^\infty a_j \phi_j \ ,
\]
where $a_j$'s are complex constants.
\end{lemma}

Proof: Since Equation (\ref{tra}) is a linear equation, the claim 
of Lemma \ref{trale} implies the current lemma. $\Box$

\subsection{The Spectrum of $L$}

Denote by $f_\tau$ the flow generated by the vector field 
($\Om_y , -\Om_x$) on $T^2$.
\begin{theorem}
Denote by $\sg_X(L)$ the spectrum of the operator $L$ (\ref{laxpair}) 
in the space $X$.
\begin{enumerate}
\item If $\Om \not \equiv \mbox{constant}$, then 
$\sg_{L^2(T^2)} (L) = iR \ $.
\item Denote by $\La$ the collection of all Lyapunov exponents of the 
flow $f_\tau$, then $\La / \{ 0 \} +iR \subset \sg_{H^1(T^2)} (L)\ $.
\end{enumerate}
\end{theorem}
For a proof of this theorem, see \cite{LS02}.

\eqnsection{3D Euler Equation}

The 3D Euler equation can be written in vorticity form,
\begin{equation}
\pa_t \Om + (u \cdot \na) \Om - (\Om \cdot \na) u = 0 \ ,
\label{3deuler}
\end{equation}
where $u = (u_1, u_2, u_3)$ is the velocity, $\Om = (\Om_1, \Om_2, \Om_3)$
is the vorticity, $\na = (\pa_x, \pa_y, \pa_z)$, 
$\Om = \na \times u$, and $\na \cdot u = 0$. $u$ can be 
represented by $\Om$ for example through Biot-Savart law.

\begin{theorem}[Childress, \cite{Chi00}]
The Lax pair of the 3D Euler equation (\ref{3deuler}) is given as
\begin{equation}
\left \{ \begin{array}{l} 
L \varphi = \la \varphi \ ,
\\
\pa_t \varphi + A \varphi = 0 \ ,
\end{array} \right.
\label{3dlaxpair}
\end{equation}
where
\[
L \varphi = (\Om \cdot \na )\varphi - (\varphi \cdot \na )\Om \ , 
\ \ \ A \varphi = (u \cdot \na )\varphi - (\varphi \cdot \na ) u \ , 
\]
$\la$ is a complex constant, and $\varphi = (\varphi_1, \varphi_2, 
\varphi_3)$ is a complex 3-vector valued function.
\end{theorem}
\begin{theorem}[Montgomery-Smith, \cite{Mon02}]
Another Lax pair of the 3D Euler equation (\ref{3deuler}) is given as
\begin{equation}
\left \{ \begin{array}{l} 
L \phi = \la \phi \ ,
\\
\pa_t \phi + A \phi = 0 \ ,
\end{array} \right.
\label{alaxpair}
\end{equation}
where
\[
L \phi = (\Om \cdot \na )\phi \ , 
\ \ \ A \varphi = (u \cdot \na )\phi \ , 
\]
$\la$ is a complex constant, and $\phi$ is a complex scalar-valued function.
\end{theorem}

\subsection{Isospectral Theory and Conservation Laws}

Denote by $H^s$ the Sobolev space $H^s(R^3)$ or $H^s(T^3)$, 
and $\| \ \|_s$ the $H^s$ norm.
\begin{theorem}
Let $\Om$ be a solution to the 3D Euler equation (\ref{3deuler}),
$\phi$ be a solution to the Lax pair (\ref{alaxpair}) at ($\Om , \la $),
then 
\[
I = \frac{\| (\Om \cdot \na )\phi  \|_s}
{\| \phi \|_s}
\]
is conservation law, i.e., I is independent of $t$. 
\end{theorem}

Proof: Take the $H^s$ norm on both side of the first equation 
in the Lax pair (\ref{3dlaxpair}), then 
\[
I = \frac{\| (\Om \cdot \na )\phi  \|_s}
{\| \phi \|_s} = |\la |.
\]
By the isospectral property of the Lax pair (\ref{3dlaxpair}), 
I is independent of $t$. $\Box$

\begin{remark}
Of course, the significance of the conservation laws comes from their 
potential in providing {\em{a priori}} bounds and establishing the global 
well-posedness of 3D Euler equation, hence, of 3D Navier-Stokes 
equation (one of the Clay problems). 
\end{remark}
\begin{lemma}
If $\{ \varphi_j \}_{j=1,2, \cdots }$ is a complete base of $H^s$,
where $\varphi_j$'s solve the Lax pair (\ref{3dlaxpair}) at different 
values of $\la$; then
\[
\Om = \sum_{j=1}^\infty a_j \varphi_j \ ,
\]
where $a_j$'s are complex constants.
\end{lemma}

Proof: The proof is completed by comparing the second equation 
in the Lax pair (\ref{3dlaxpair}) and the 3D Euler equation. $\Box$

\bibliography{Iso}

\begin{thebibliography}{1}

\bibitem{Li01a}
Y.~Li.
\newblock A {L}ax {P}air for the {T}wo {D}imensional {E}uler {E}quation.
\newblock {\em J. Math. Phys.}, 42, No.8:3552--3553, 2001.

\bibitem{LY00}
Y.~Li and A.~Yurov.
\newblock Lax {P}airs and {D}arboux {T}ransformations for {E}uler {E}quations.
\newblock {\em Submitted, available at:
  http://xxx.lanl.gov/abs/math.AP/0101214, or
  http://www.math.missouri.edu/\~{}cli}, 2000.

\bibitem{LS02}
Y.~Latushkin and R.~Shvidkoy.
\newblock Spectral {P}roperties of {L}inearized 2{D} {E}uler {E}quation and
  {E}volution {S}emigroup.
\newblock {\em Preprint (available upon request)}, 2002.

\bibitem{Chi00}
S.~Childress.
\newblock {\em personal communication}, 2000.

\bibitem{Mon02}
S.~Montgomery-Smith.
\newblock {\em personal communication}, 2002.

\end{thebibliography}

\end{document}